\documentclass[12pt,a4paper,twoside,final,notitlepage, reqno]{article}

\usepackage[english]{babel}
\usepackage[latin1]{inputenc}
\usepackage[a4paper,left=3.5cm,right=3.5cm]{geometry}
\usepackage{amssymb}
\usepackage{amsthm} 
\usepackage{amsmath}
\usepackage{mathtools}
\usepackage{amscd}
\usepackage{geometry}
\usepackage{graphics,graphicx}
\usepackage{epstopdf}
\usepackage[usenames, dvipsnames]{color}
\usepackage[textsize=small]{todonotes}
\usepackage{booktabs}
\usepackage{subfigure} 
\usepackage{siunitx}

\graphicspath{{figs/}}
\makeatletter
\def\input@path{{figs/}}
\makeatother

\theoremstyle{definition}
\theoremstyle{definition}
\theoremstyle{definition}
\usepackage{esdiff}
\renewcommand{\d}[1]{\, \mathrm{d} #1}
\def \dx {\d{x}}
\def \ds {\d{s}}

\newcommand{\R}{\mathbb{R}}

\renewcommand{\phi}{\varphi}

\def\sizeg{  \footnotesize }

\usepackage{authblk}

\title{\bf{\Large{
      Approximation of the Ventcel problem,
      \\[5pt]
      numerical results
    }}
}

\author[1]{Marc Dambrine \thanks{marc.dambrine@univ-pau.fr}}
\author[1]{Charles Pierre \thanks{charles.pierre@univ-pau.fr}}
\affil[1]{
  Laboratoire  de Math\'ematiques et de leurs Applications, 
  , UMR CNRS 5142, \protect \\
  Universit\'e de Pau et des Pays de l'Adour, France.
}

\usepackage{fancyhdr}

\begin{document} 

\date{December,  2017}

\maketitle
% \noindent
% \keywords{}
% \\ \\ 
% \subjclass{}
% \\ \\
% \acknow{}
% \\ \\
%
%
\begin{abstract}
Report on the numerical approximation of the Ventcel problem. The Ventcel problem is a 3D eigenvalue problem involving a surface differential operator on the domain boundary: the Laplace-Beltrami operator.
\\
We present in the first section the problem statement together with its finite element approximation, the code machinery used for its resolution is also presented here.
\\
The last section presents the obtained numerical results. These results are quite unexpected for us. Either super-converging for $P^1$ Lagrange finite elements or under converging for $P^2$ and $P^3$.
\\
The remaining sections 2 and 3 provide numerical results either for the classical Laplace  or for the Laplace-Beltrami operator numerical approximation. These examples being aimed to validate the code implementation.
\end{abstract}

\section{Introduction}

We present in this report an algorithm for solving the Ventcel problem \cite{ventcel-56,ventcel-59}. The code is detailed and validated considering various classical test cases. The resolution of the Ventcel problem on the unit ball is used as a test case to provide  a convergence numerical analysis of the method.

\subsection{The Ventcel problem}
\label{sec:ventcel-pb}
The Ventcel problem is the following eigenvalue problem. Let $\Omega$ denote some bounded smooth domain in $\R^3$. We search for the eigenvalues $\lambda$ and for the associated eigenfunctions $u$ satisfying,
\begin{equation}
  \label{eq:ventcel}
  \Delta u = 0 \quad \text{on}\quad \Omega,\quad \Delta_B u - \partial_n u + \lambda u =0\quad \text{on}\quad \partial\Omega.
\end{equation}
Here $\mathbf{n}$ denotes the outward unit normal to $\partial\Omega$ and so $\partial_n u$  the derivative of $u$ on the normal direction to the boundary and $\Delta_B u$ stands for the surface Laplacian of $u$ on the boundary $\partial\Omega$, i.e. the Laplace-Beltrami operator.

This problem has a simple weak formulation \eqref{eq:ventcel-weak}. Multiplying $\Delta u$ by a test function $v$ and integrating over $\Omega$ and assuming sufficient regularity we get by integration by part,
\begin{displaymath}
  \int_\Omega \nabla u\cdot\nabla v \dx - \int_{\partial\Omega} \partial_n u~ v \ds =0.
\end{displaymath}
By substituting the boundary condition in \eqref{eq:ventcel},
\begin{displaymath}
  \int_\Omega \nabla u\cdot\nabla v \dx -
  \int_{\partial\Omega} \left(  \Delta_B u + \lambda u\right)v\ds.
\end{displaymath}
We obtain by a second integration by part on the boundary $\partial\Omega$,
\begin{equation}
  \label{eq:ventcel-weak}
  \int_\Omega \nabla u\cdot\nabla v \dx + \int_{\partial\Omega}  \nabla_T u\cdot\nabla_T v \ds
  = \lambda \int_{\partial\Omega} u v \ds,
\end{equation}
with $\nabla _T u$ the tangential gradient of $u$ on $\partial\Omega$, i.e. $\nabla _T u= \nabla u - \partial_n u ~\mathbf{n}$.

On the Hilbert space $V=\{u\in {\rm H}^1(\Omega),~\nabla_T u \in {\rm L}^2(\partial\Omega)\}$ the classical theory applies as exposed in \cite{dambrine-kateb-2012} providing an orthogonal set of eigenfunctions associated with the eigenvalues $0=\lambda_0 < \lambda_1 \le \dots$. The first eigenvalues $\lambda_0=0$ is associated with the (one dimensional) eigenspace of the constants. We will focus on the approximation of $\lambda_1$.

The following discretisation of \eqref{eq:ventcel-weak} will be used. Let $\mathcal{M}$ be some tetrahedral mesh of $\Omega$. 
We will denote by $\Omega_h$ the computation domain (made of all tetrahedra in the mesh $\mathcal{M}$).
We consider the classical Lagrange finite element spaces $V_h=P^k(\mathcal{M})$. We search for $u\in V_h$ and $\lambda\in\R$ so that, for all $v\in V_h$ we have:
\begin{equation}
  \label{eq:ventcel-disc-0}
  \int_{\Omega_h} \nabla u\cdot\nabla v \dx + \int_{\partial\Omega_h}  \nabla_T u\cdot\nabla_T v \ds
  = \lambda \int_{\partial\Omega_h} u_h v \ds.
\end{equation}
Considering a canonical bases of $V_h$ we identify $u_h$ with its vectorial representation $U\in\R^N$, with $N$ the number of degrees of freedom (the dimension of $V_h$). We then consider the matrices $S_3$, $S_2$ and $M_2$ (mass and stiffness matrices) representing the products on $V_h\times V_h$ in \eqref{eq:ventcel-disc-0}:
\begin{align}
  \label{eq:prod-S3}
  (u,v)\mapsto &
  \int_{\Omega_h} \nabla u\cdot\nabla v \dx
    = U^T S_3 V,
  \\
  \label{eq:prod-S2}
  (u,v)\mapsto  &
  \int_{\partial\Omega_h}  \nabla_T u\cdot\nabla_T v \ds
  = U^T S_2 V,
  \\
  \label{eq:prod-M2}
  (u,v)\mapsto  &
  \int_{\partial\Omega_h} u v \ds
  = U^T M_2 V.
\end{align}

Problem \eqref{eq:ventcel-disc-0} under matricial form is the following generalized eigenvalue problem,
\begin{equation}
  \label{eq:ventcel-disc}
  (S_3+S_2)U =\lambda M_2 U.
\end{equation}

 The matrix $S_3+S_2$ is symmetric positive semi definite (with a one dimensional kernel made of the constant vectors $U=c$ representing the constant functions). The matrix $M_2$ also is symmetric positive semi definite but with a high dimensional kernel made of all functions in $V_h$ vanishing on the boundary $\partial\Omega_h$.

\subsection{Implementation}
\label{sec:impl}
A \textit{python} code is used for the numerical resolution  of \eqref{eq:ventcel-disc}. The assembling of the finite element matrices is made using \textit{Getfem++}
\footnote{Getfem++: an open-source finite element library, http://download.gna.org/getfem/html/homepage/index.html}.
\\
The iterative Lanczos method (see e.g. Saad \cite{saad-vp}) is used to solve the eigenvalue problem. 
Because we are interested in the smallest eigenvalues and since the matrix $S_3+S_2$ is only semi definite, we consider the shifted invert variant of this method.
Actually the shifted matrix $S_3+S_2-\sigma M_2$ is symmetric positive definite for any negative value of the shift parameter $\sigma$, fixed to  $\sigma=-1$ here.
In practice the ARPACK library\footnote{ARPACK, Arnoldi Package, http://www.caam.rice.edu/software/ARPACK/} has been used for this. 
\\
Each Lanczos algorithm iteration requires two kind of advanced operations:
solving the linear system $(S_3+S_2- M_2) X = Y$ (shifted system matrix) and performing matrix-vector multiplication $X\mapsto  M_2 X$. These two operations are executed using the \textit{PETSc} library\footnote{PETSc, Portable, Extensible Toolkit for Scientific Computation, http://www.mcs.anl.gov/petsc/}.
\\
The largest amount of CPU consumption is devoted to the linear system resolution. A conjugate gradient algorithm is used with an incomplete Cholesky (fill in 3) preconditioning, see e.g. Saad \cite{saad} for precisions. The system is solved with a tolerance of $10^{-11}$.
\\
The residual at end of the Lanczos algorithm is asked to be smaller than $10^{-10}$. This residual is algebraic, we additionally re-computed an ${\rm L}^2$ residual: denoting $\lambda_i$ the $i^{th}$ eigenvalue and $U_i$ the associated eigenvector, the ${\rm L}^2$ residual is defined as,
\begin{displaymath}
  \dfrac{  \Vert (S_3+S_2)U_i - \lambda_i M_2 U_i\Vert_{{\rm L}^2}}
  {   \Vert  U_i\Vert_{{\rm L}^2}},
\end{displaymath}
assuming the identification $V_h \simeq \R^N$.
This functional residual is also controlled to remain below $10^{-9}$.
\\
Eventually all meshes have been built with the software \textit{GMsh}\footnote{Gmsh: a three-dimensional finite element mesh generator, http://geuz.org/gmsh/}.
\section{Code validation}

We consider here the resolution of several classical problems aimed to validate the code implementation: matrix assembling, eigenproblem solver and error analysis. On each example our purpose is to recover the correct method's convergence rate as predicted by classical theories.
\subsection{Matrix assembling validation}
\label{sec:mat-assemb-valid}
The three matrices $S_3$, $S_2$ and $M_2$ associated with the products (\ref{eq:prod-S3}), (\ref{eq:prod-S2}) and (\ref{eq:prod-M2}) need to be computed.
The matrix $S_3$ is a classical stiffness matrix on the three dimensional domain $\Omega_h$. 
Conversely the matrices $S_2$ and $M_2$ are not as classical. They correspond to two dimensional stiffness and mass matrices but associated with the non flat domain $\partial\Omega_h$. 
\\
Remark nevertheless that the assembling of $S_2$ (and similarly for $M_2$) consist of a loop over all element faces $E$ of the mesh $\mathcal{M}$ included in $\partial\Omega_h$. Precisely the elements of $\mathcal{M}$ are tetrahedra, their faces $E$ thus are triangles and these triangles subsets of $\partial\Omega_h$ form a partition of $\partial\Omega_h$. On every such triangle $E$ is assembled the local matrix corresponding to the product,
\begin{displaymath}
  (u,v)\mapsto \int_E \nabla _T u\cdot\nabla _T v \ds = \int_E \nabla u_{|E}\cdot\nabla v_{|E} \ds.
\end{displaymath}
On one hand $E$ is an affine deformation of the reference triangle in dimension 2 and on the other hand $u_{|E}, v_{|E}\in P^k(E)$. Therefore these local matrices indeed are local matrices for a classical 2 dimensional stiffness matrix. As a result the assembling of $S_2$ and similarly of $M_2$ do not need particular software but only a 2D finite element library.

Taking the previous commentary into account, the validation for the assembling of $S_3$, $S_2$ and $M_2$ can be performed on flat domains $\omega\subset \R^d$, $d=2,3$. We fix $\omega=(0,1)^d$ and solve the classical elliptic problem,
\begin{displaymath}
  \label{eq:lap-pb}
  -\Delta u + u = f\quad \text{on}\quad \omega \quad \text{and}\quad\partial_n u = 0\quad \text{on}\quad \partial\omega,
\end{displaymath}
for the right hand side $f=(2\pi^2+1)\cos(\pi x) \cos(\pi y)$ and $f=(3\pi^2+1)\cos(\pi x)  \cos(\pi y)\cos(\pi z)$ in dimension 2 and 3 respectively. The associated exact solutions are $u=\cos(\pi x) \cos(\pi y)$ and $u=\cos(\pi x) \cos(\pi y)\cos(\pi z)$ respectively. The numerical solution $u_h$ is given by,
\begin{displaymath}
  S_d U_h + M_d U_h = M_d F~,\quad d=2,3.
\end{displaymath}
We analyze the relative errors between $u$ and $u_h$ both in ${\rm L}^2$-norm and ${\rm H}^1$-semi norm,
\begin{displaymath}
  e_{_{{\rm L}^2}}=\dfrac{\Vert u-u_h\Vert _{{\rm L}^2(\omega)}}{\Vert u\Vert _{{\rm L}^2(\omega)}}~,
  \quad 
  e_{_{{\rm H}^1}}=\dfrac{\Vert u-u_h\Vert _{{\rm H}^1(\omega)}}{\Vert u\Vert _{{\rm H}^1(\omega)}}.
\end{displaymath}
considering a series of refined meshes and using the $P^k(\omega)$ finite element method for $k=1,2,3$, we recovered a $k$ (resp. $k+1$) order of convergence in ${\rm H}^1$ (resp. ${\rm L}^2$) norm as presented in table \ref{tab:res-lap}. These results are in complete agreement with the classical theory (see e.g. Ciarlet \cite{ciarlet-ef-elliptic}) and this test fully validates a correct assembling of the desired matrices.
\begin{table}[h!]
  \centering
  \begin{tabular}{ccc}

  \begin{tabular}{|r|c|c|c|}
    \hline
    $d=2\quad \quad \quad $ &$~~P^1~~$&$~~P^2~~$&$~~P^3~~$
    \\ \hline
    $e_{_{{\rm H}^1}}$ &1.0&2.1&3.0
    \\ \hline
    $e_{_{{\rm L}^2}}$ &2.0&3.1&4.0
    \\ \hline
  \end{tabular}
  &  ~~~~~~&
  \begin{tabular}{|r|c|c|c|}
    \hline
    $d=3\quad \quad \quad $ &$~~P^1~~$&$~~P^2~~$&$~~P^3~~$
    \\ \hline
    $e_{_{{\rm H}^1}}$ &1.0&1.9&2.8
    \\ \hline
    $e_{_{{\rm L}^2}}$ &1.9&3.0&4.0
    \\ \hline    
  \end{tabular}   
  \end{tabular}
  \caption{Computed orders of convergence for problem \eqref{eq:lap-pb}}
  \label{tab:res-lap}
\end{table}

\subsection{Eigenvalue problem solver validation}
\label{sec:eigen-solver-valid}
We test the eigenvalue problem solver described in section \ref{sec:impl} considering the Laplace eigenproblem,
\begin{equation}
  \label{eq:lap-vp}
  -\Delta u = \lambda u\quad \text{on}\quad \omega \quad \text{and}\quad \partial_n u = 0 \quad \text{on}\quad \partial\omega,
\end{equation}
on the same square or cubic geometry $\omega$ as in the previous section. 
\\
The first eigenvalue is $\lambda_0=0$ with eigenspace the constant functions. The first non-zero eigenvalue is $\lambda_1=\pi^2$ with multiplicity $d$ and with eigenspace $E_1=\text{Span}\left(\cos(\pi x), \cos(\pi y)\right)$ or $E_1=\text{Span}\left(\cos(\pi x), \cos(\pi y), \cos(\pi z)\right)$ for $d=2$ or $d=3$.

According to the dimension $d$ the numerical approximation for \eqref{eq:lap-vp} is,
\begin{displaymath}
  S_d U_h = \lambda_h M_d U_h.
\end{displaymath}
The first non-zero numerical eigenvalue $\lambda_{h,1}$ is computed to analyze the relative error,
\begin{displaymath}
  e_\lambda = \dfrac{|\lambda_{h,1}-\lambda_{1}|}{\lambda_{1}}.
\end{displaymath}
One associated eigenvector $U_{h,1}$ of ${\rm L}^2$-norm equal to 1 (for normalization) is considered to compute the following ${\rm H}^1$ and ${\rm L}^2$ errors,
\begin{displaymath}
  e_{_{{\rm L}^2}}=\Vert u_{h.1}- pu_{h,1}\Vert _{{\rm L}^2(\omega)}
  \quad 
  e_{_{{\rm H}^1}}=\Vert u_{h,1} - pu_{h,1}\Vert _{{\rm H}^1(\omega)},
\end{displaymath}
where $p$ is the ${\rm L}^2$-orthogonal projection  onto the eigenspace $E_1$. 

\begin{table}[h!]
  \centering
  \begin{tabular}{ccc}
    
    \begin{tabular}{|r|c|c|c|}
      \hline
      $d=2\quad \quad \quad $ &$~~P^1~~$&$~~P^2~~$&$~~P^3~~$
      \\ \hline
      $e_{_{{\rm H}^1}}$ &1.0&2.0&3.1
      \\ \hline
      $e_{_{{\rm L}^2}}$ &2.0&3.0&4.1
      \\ \hline
      $e_{\lambda}$ &2.0&4.2&5.9
      \\ \hline
    \end{tabular}
    &  ~~~~~~&
    \begin{tabular}{|r|c|c|c|}
      \hline
      $d=3\quad \quad \quad $ &$~~P^1~~$&$~~P^2~~$&$~~P^3~~$
      \\ \hline
      $e_{_{{\rm H}^1}}$ &1.0&2.0&2.7
      \\ \hline
      $e_{_{{\rm L}^2}}$ &1.9&3.0&3.8
      \\ \hline
      $e_{\lambda}$ &1.95&4.0&5.9
      \\ \hline    
    \end{tabular}   
  \end{tabular}
  \caption{Computed orders of convergence for problem \eqref{eq:lap-vp}}
  \label{tab:res-lap-vp}
\end{table}

The order of convergence is computed considering a series of refined meshes and several Lagrange finite element spaces $P^k(\omega)$, they are reported in table \ref{tab:res-lap-vp}.
They are in full agreement with the theoretical orders presented e.g. in Babu{\v{s}}hka and Osborn \cite{babuska-osborne-eigen}~: using $P^k$ finite elements provides a convergence of order $k$ for $e_{_{{\rm H}^1}}$, $k+1$ for $e_{_{{\rm L}^2}}$ and $2k$ for $e_\lambda$.

\section{Curved surface effects}

In this section we illustrate the influence of computing the two dimensional stiffness and mass matrices $S_2$ and $M_2$. We consider again the geometrical situation of the Ventcel problem in section \ref{sec:ventcel-pb}: $\Omega$ is an open smooth domain in $\R^3$ with boundary $\Gamma$. $\mathcal{M}$ is a tetrahedral mesh of $\Omega$, $\Omega_h$ is the  domain of the mesh $\mathcal{M}$, i.e. the union of its element. The matrices $S_2$ and $M_2$ are computed on  $\Gamma_h=\partial\Omega_h$. Note that $\Gamma_h$ itself is the  domain of a triangular mesh $\mathcal{M}_2$ of $\Gamma$, the element of $\mathcal{M}_2$ are the boundary faces of $\mathcal{M}$.
\\
As developed in section \ref{sec:mat-assemb-valid}, the assembling of $S_2$ and $M_2$ on the curved domain $\Gamma_h$ (i.e. embedded in $\R^3$) is essentially the same as their assembling in a flat domain (a two dimensional domain in $\R^2$). We however check here that the third dimension $z$ is taken into account for this assembling.
\subsection{Laplace-Beltrami problem}
We will consider the elliptic problem
\begin{equation}
  \label{eq:lap-bel-pb}
  -\Delta_B u + u = f~\quad \text{on}\quad \Gamma,
\end{equation}
no boundary condition are needed here ($\Gamma$ has no boundary!). 
This problem is well posed : existence and uniqueness of a solution $u\in {\rm H}^1(\Gamma)$ for any $f\in {\rm L}^2(\Gamma)$.
\\
The finite element approximation of \eqref{eq:lap-bel-pb} has been theorized by Demlow in \cite{Demlow}.
We need some notations. To $x\in\Gamma$ is associated its unit outer normal $\mathbf{n}(x)$. We consider a tubular neighborhood $\omega$ of $\Gamma$ so that for all $x\in\omega$, there exists a unique $p(x)\in\Gamma$ so that $x-p(x) = \lambda(x)\mathbf{n}(p(x))$ ($p(x)$ is an orthogonal projection of $x$ on $\Gamma$) and so that the segment $[x,p(x)]\subset \omega$.
A function $f:~\Gamma\longrightarrow \R$ can be extended to a function $f^e:~\omega\longrightarrow \R$ by: $f^e(x)=f(p(x))$. We assume that $\Gamma_h\subset  \omega$. A function $f:~\Gamma\longrightarrow \R$ can then be lifted to a function $f^l:~\Gamma_h\longrightarrow \R$ by $f^l=f^e_{|\Gamma_h}$.
\\
The lift operation allows to compare the exact solution $u$ defined on $\Gamma$ with a numerical approximation $u_h$ defined on $\Gamma_h$. Note that it would have also been possible to lift $u_h$ to a function $u_h^l$ on $\Gamma$ by $u_h^l(p(x))=u_h(x)$. Demlow showed in \cite{Demlow} that analyzing the error in terms of $u^l-u_h$ (lift of $u$ to $\Gamma_h$) or in terms of $u - u_h^l$ (lift of $u_h$ to $\Gamma$) are equivalent, we choose the first strategy for its practical convenience.

Equation \eqref{eq:lap-bel-pb} is discretised on $V_h=P^k(\mathcal{M}_2)$ as,
\begin{equation}
  \label{eq:lap-bel-pb-h}
  (S_2+M_2)U_h = M_2 F^l.
\end{equation}
The ${\rm L}^2$ and ${\rm H}^1$ approximation  errors are alternatively defined as,
\begin{displaymath}
  e_{_{{\rm L}^2}}=\Vert u^l-u_h\Vert _{{\rm L}^2(\Gamma_h)}~,
  \quad 
  e_{_{{\rm H}^1}}=\Vert u^l-u_h\Vert _{{\rm H}^1(\Gamma_h)}.
\end{displaymath}
The convergence properties of this scheme are quite different from those on flat domains illustrated in section \ref{sec:mat-assemb-valid}.
The reason for this is analyzed in Demlow \cite{Demlow}: the curved surface $\Gamma$ is approximated by the surface $\Gamma_h$ that is the boundary of a polyhedral and thus piecewise flat (made of triangles). Approximating a curved surface by a (piecewise) linear one induces upper bound errors for \eqref{eq:lap-bel-pb},
\begin{displaymath}
  e_{_{{\rm H}^1}} = O(h^k + h^2)~,\quad e_{_{{\rm L}^2}} = O(h^{k+1} + h^2),
\end{displaymath}
for $P^k$ finite elements. Therefore a saturation of the convergence order to 2 is predicted.

The numerical scheme \eqref{eq:lap-bel-pb-h} has been implemented for the sphere. In this case the projection $p$ is very simple, $p(x)=x/||x||$ on $\omega=\R^3-\{0\}$. The right hand side to \eqref{eq:lap-bel-pb} is set to $f=x(2+x)\exp(x)$ so that the exact solution is $u=\exp(x)$. The results are reported in table \ref{tab:res-lap-bel}. The saturation of the convergence order to $h^2$ is clearly seen, in agreement with Demlow convergence analysis \cite{Demlow}.
\begin{table}[h!]
  \centering
  \begin{tabular}{ccc}
    
    \begin{tabular}{|r|c|c|c|}
      \hline
      &$~~P^1~~$&$~~P^2~~$&$~~P^3~~$
      \\ \hline
      $e_{_{{\rm H}^1}}$ &1.1&2.0&2.0
      \\ \hline
      $e_{_{{\rm L}^2}}$ &2.0&2.0&2.0
      \\ \hline
    \end{tabular}
    &  ~~~~~~~~~~~~~&
    \begin{tabular}{|r|c|c|c|}
      \hline
      &$~~P^1~~$&$~~P^2~~$&$~~P^3~~$
      \\ \hline
      $e_{_{{\rm H}^1}}$ &1&2&2
      \\ \hline
      $e_{_{{\rm L}^2}}$ &2&  2  &  2
      \\ \hline
      $e_{\lambda}$ &2&2&2
      \\ \hline    
    \end{tabular}   
  \end{tabular}
  \caption{Computed orders of convergence for the Laplace-Beltrami problem \eqref{eq:lap-bel-pb} (left) and for the Laplace-Beltrami eigenvalue problem \eqref{eq:lap-bel-vp}.}
  \label{tab:res-lap-bel}
\end{table}

\subsection{Laplace-Beltrami eigenvalue problem}

The Laplace-Beltrami eigenproblem is considered on the sphere,
\begin{equation}
  \label{eq:lap-bel-vp}
  -\Delta_B u =\lambda u \quad \text{on}\quad \Gamma.
\end{equation}
Its discretisation takes the form: find $U_h\in V_h=P^k(\Gamma_h)$ and $\lambda_h\in \R$ so that,
\begin{displaymath}
  S_2 U_h = \lambda_h M_2 U_h.
\end{displaymath}
Problem \eqref{eq:lap-bel-vp} has $\lambda_0=0$ for eigenvalue associated to the eigenspace of constant functions. The first non-zero eigenvalue is $\lambda_1=2$ of multiplicity 3 with eigenspace $E_1=\text{Span}(x,y,z)$ the restriction of the linear functions in $\R^3$ to the sphere.

The space $E_1$ is lifted to $\Gamma_h$ as in the previous section. The lifted space $E_i^l$ is the vector space of functions of the form $X=(x,y,z) \mapsto (\alpha x  + \beta y+\delta z)/||X||$ with $||X||=(x^2+y^2+z^2)^{1/2}$. The orthogonal projector $p$ from $V_h$ onto $E_1^l$ is considered to define the numerical errors,
\begin{displaymath}
  e_\lambda = \dfrac{|\lambda_{h,1}-\lambda_{1}|}{\lambda_{1}},
  \quad 
  e_{_{{\rm L}^2}}=\dfrac{
    \Vert u_{h,1}- pu_{h,1}\Vert _{{\rm L}^2(\Gamma_h)}}{\Vert u_{h,1}\Vert _{{\rm L}^2(\Gamma_h)}}
  \quad \text{and}\quad 
  e_{_{{\rm H}^1}}=\dfrac{\Vert u_{h,1} - pu_{h,1}\Vert _{{\rm H}^1(\Gamma_h)}}{\Vert u_{h,1}\Vert _{{\rm L}^2(\Gamma_h)}},
\end{displaymath}
on the first computed non zero eigenvalue $\lambda_{h,1}$ and associated eigenfunction $u_{h,1}/\Vert u_{h,1}\Vert _{{\rm L}^2(\Gamma_h)}$ normalized in ${\rm L}^2$-norm.

The results are reported in table \ref{tab:res-lap-bel} on the right. The convergence in ${\rm H}^1$-norm is in $O(h^k+h^2)$ for $P^k$ Lagrange finite element which is consistent with the previous subsection assertions.
Order 2 convergence in ${\rm L}^2$-norm is observed in all cases $P^k$, $k=1,~3$.
\\
The eigenvalue convergence clearly shows a saturation to order 2 convergence. 
In  section \ref{sec:eigen-solver-valid} we quote the analysis of  Babu{\v{s}}ka and Osborn \cite{babuska-osborne-eigen} who showed an order $2k$ convergence for the eigenvalues with $k$ the convergence order of the eigenfunctions in ${\rm H}^1$-norm. This analysis does not apply here since it is restricted to Galerkin type approximations with $V_h\subset  V$. In the present case this is no longer true because $V$ and $V_h$ are function spaces associated  to different domains, $\Gamma$ and $\Gamma_h$ respectively.

{\bf Conclusion.}
The numerical results showed in this section together with the theoretical analysis of Demlow \cite{Demlow} show that one cannot get better than an order 2 convergence when considering a piecewise affine mesh $\Gamma_h$ of the curved boundary $\Gamma$ of the domain $\Omega$. This point also was addressed in \cite{Demlow} where piecewise polynomial interpolation of order $p$ $\Gamma^p_h$  of $\Gamma$ are considered providing now a saturation of the convergence order to $p+1$.

\section{The Ventcel problem numerical convergence}
\begin{figure}[h!]
  \centering
  \begin{tabular}{c}
    \hspace{-2.6cm}
    \input{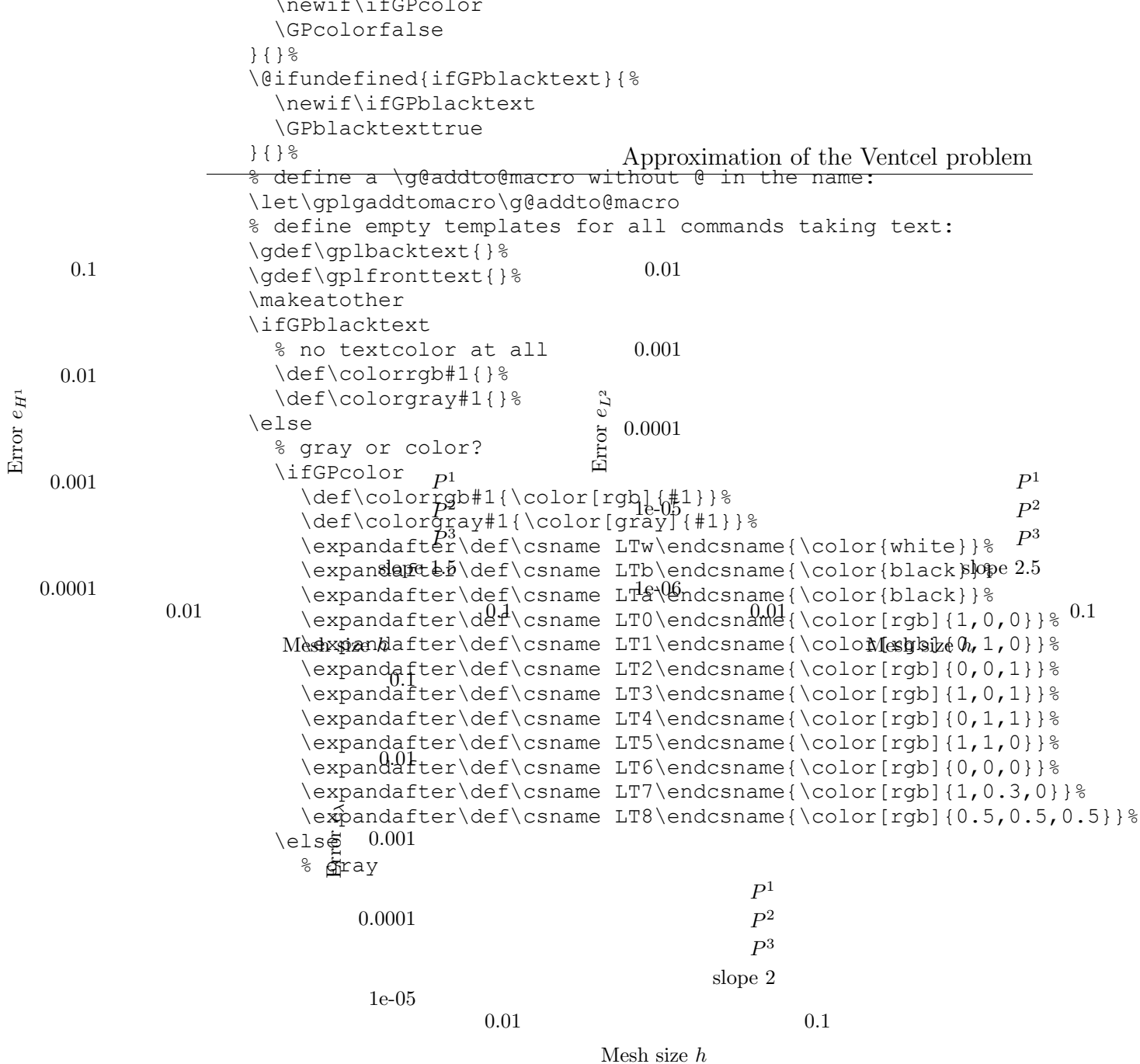}
  \end{tabular}\hspace{-2cm}
  \caption{Convergence for the Ventcel problem in logarithmic scale. Above, errors on the first computed eigenfunction $u_{h,1}$ in ${\rm H}^1$ (left) and ${\rm L}^2$ norms. Below, errors on the first computed eigenvalue $\lambda_1$.}
  \label{fig:ventcel-vp1}
\end{figure}
\begin{table}[h!]
  \centering
  \begin{tabular}{c}
    
    \begin{tabular}{|r|c|c|c|}
      \hline
      &$~~P^1~~$&$~~P^2~~$&$~~P^3~~$
      \\ \hline
      $e_{_{{\rm H}^1}}$ &1.&1.5&1.5
      \\ \hline
      $e_{_{{\rm L}^2}}$ &2.&2.5&2.5
      \\ \hline
      $e_{\lambda}$ &2.&2.&2.
      \\ \hline
    \end{tabular}
  \end{tabular}
  \caption{Convergence orders for the Ventcel problem}
  \label{tab:res-ventcel}
\end{table}
We analyze here the numerical approximation of the Ventcel problem exposed in section~ \ref{sec:ventcel-pb} using the numerical scheme (\ref{eq:ventcel-disc}). 

The  domain $\Omega$ is set to the unit ball. 
In this particular case the eigenvalues/eigenfunctions have been computed in \cite{dambrine-kateb-lamboley}. 
The $n^{th}$ eigenvalue $\lambda_n$ for $n\ge 0$ is given by
$  \lambda_n = n^2+2n$, with multiplicity $2n+1$ and with eigenspace $E_n$ the space of harmonic functions of order $n$. 
\\
The numerical domain $\Omega_h$ here is a subset (by convexity)  of $\Omega$. 
Therefore functions $f$ on $\Omega$ simply are lifted to functions $f^l$ on $\Omega_h$ by restriction: $f^l = f_{|\Omega_h}$. 
\\
We consider a numerical approximation $\lambda_{n,h}$ of the $n^{th}$ non zero eigenvalue $\lambda_n$ (since $\lambda_n$ is of multiplicity $2n+1$, we numerically get $2n+1$ approximations of $\lambda_n$, we fix one of them). We denote $U_{n,h}$ an associated eigenfunction, so that $\Vert U_{n,h}\Vert_{{\rm L}^2(\Omega_h)}=1$ (${\rm L}^2$ normalisation).
We introduce the orthogonal projection (for the ${\rm L}^2$-scalar product on $\Omega_h$) $p~:V_h\longrightarrow E^l_n$, where here $E_n^l$ is the space made of the restrictions of all functions in $E_n$ to $\Omega_h$. 
\\
with these notations we define the numerical errors,
\begin{displaymath}
  e_{\lambda_n} = \dfrac{|\lambda_{n,h}-\lambda_{n}|}{\lambda_{n}},
  \quad 
  e_{_{{\rm L}^2}}=
    \Vert U_{h,h}- pU_{n,h}\Vert _{{\rm L}^2(\Omega_h)}
  \quad \text{and}\quad 
  e_{_{{\rm H}^1}}=\Vert \nabla U_{n,h} - \nabla  pU_{n,h}\Vert _{{\rm L}^2(\Omega_h)}
\end{displaymath}
on the $n^{th}$ computed non zero eigenvalue $\lambda_{h,n}$ and associated eigenfunction $u_{h,n}/\Vert u_{h,n}\Vert _{{\rm L}^2(\Omega_h)}$ normalized in ${\rm L}^2$-norm.

Numerical results are reported in table \ref{tab:res-ventcel} for $V_h=P^k(\Omega)$ and $k=1$, 2 and 3. 
 \\
 The convergence towards the first non-zero eigenvalue $l_1$ displays the expected order 2 of convergence. 
All other results are quite unexpected. 
Firstly for the convergence towards the first eigenfunction $u_{1}$. 
 The $P^1$ scheme is over-converging in ${\rm H}^1$-norm (order 1.5 whereas order 1 was expected). 
Meanwhile the $P^2$ and $P^3$ schemes are under-converging in ${\rm H}^1$-norm (order 1.5 whereas order 2 was expected). 
The three schemes are over-converging in ${\rm L}^2$-norm, displaying a convergence order of 2.5 when order 2  was expected.
 \\
The error behaviors with respect to the mesh size are depicted on figure \ref{fig:ventcel-vp1}. On top of the previous remarks, one supplementary abnormality is expressed here. On the three plots the $P^1$ scheme is the most accurate one. The $P^3$ scheme being moreover by far the worst for the ${\rm L}^2$ error $e_{_{{\rm L}^2}}$ and for the error on the eigenvalue.

% \begin{figure}[h!]
%   \centering
%   \begin{tabular}{c}
%     \hspace{-2.6cm}
%     \input{conv_Ventcel_vp2.tex}
%   \end{tabular}\hspace{-2cm}
%   \caption{Same thing as figure \ref{fig:ventcel-vp1} for the second non zero eigenvalue/eigen function (with multiplicity 5)}
%   \label{fig:ventcel-vp2}
% \end{figure}

% \begin{figure}[h!]
%   \centering
%   \begin{tabular}{c}
%     \hspace{-2.6cm}
%     \input{conv_Ventcel_vp3.tex}
%   \end{tabular}\hspace{-2cm}
%   \caption{Same thing as figure \ref{fig:ventcel-vp1} for the third non zero eigenvalue/eigen function (with multiplicity 7)}
%   \label{fig:ventcel-vp3}
% \end{figure}

\begin{figure}[h!]
  \centering
  \begin{tabular}{c}
    \hspace{-2.6cm}
    \input{conv_Ventcel_vp4.txt}
  \end{tabular}\hspace{-2cm}
  \caption{Same thing as figure \ref{fig:ventcel-vp1} for the fourth non zero eigenvalue/eigen function (with multiplicity 9)}
  \label{fig:ventcel-vp4}
\end{figure}

\begin{figure}[h!]
  \centering
  \begin{tabular}{c}
    \hspace{-2.6cm}
    \input{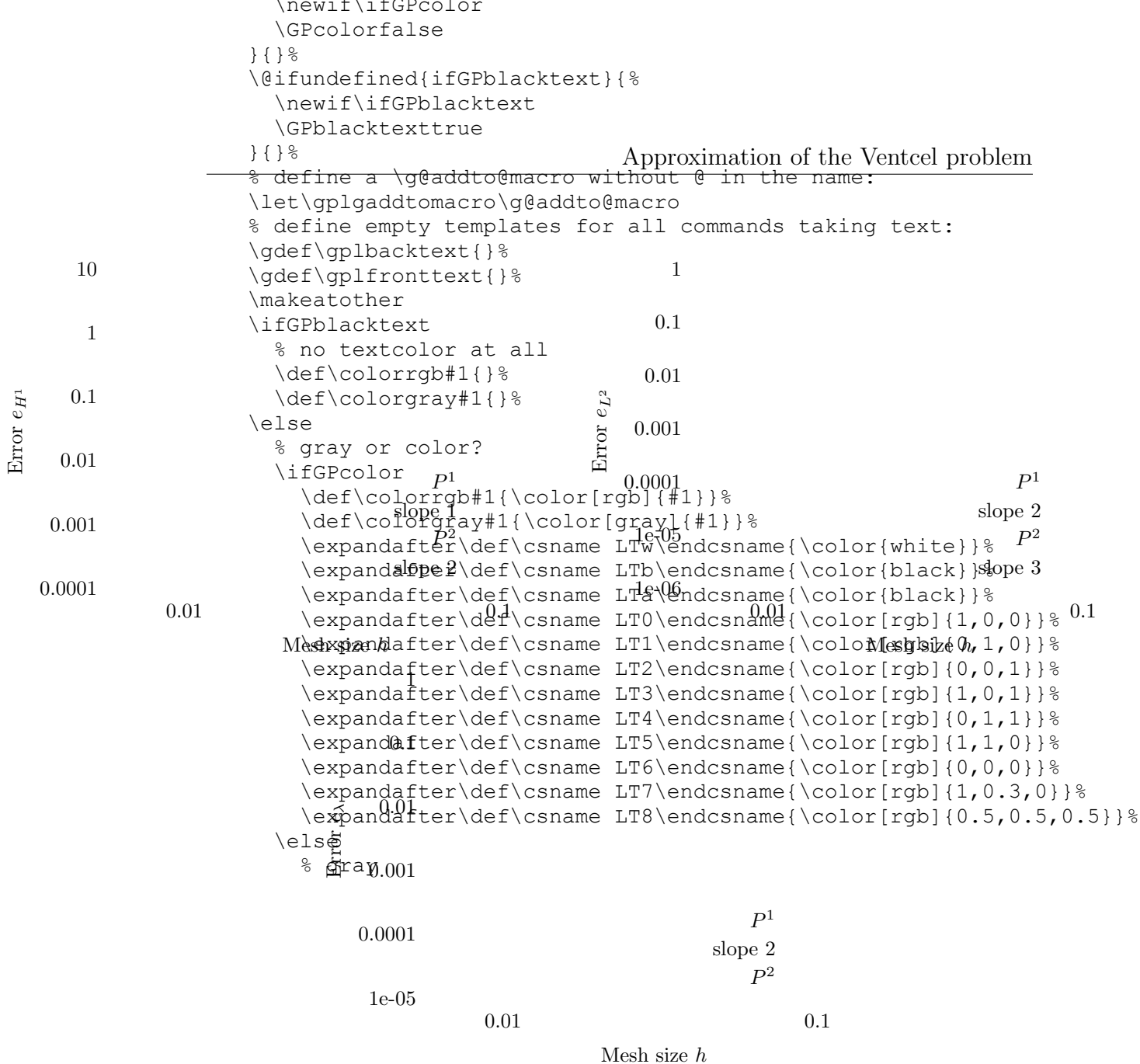}
  \end{tabular}\hspace{-2cm}
  \caption{Same thing as figure \ref{fig:ventcel-vp1} for the modified Ventcell problem and for for the fourth non zero eigenvalue/eigen function (with multiplicity 9)}
  \label{fig:ventcel-mod1-vp4}
\end{figure}

\begin{figure}[h!]
  \centering
  \begin{tabular}{c}
    \hspace{-2.6cm}
    \input{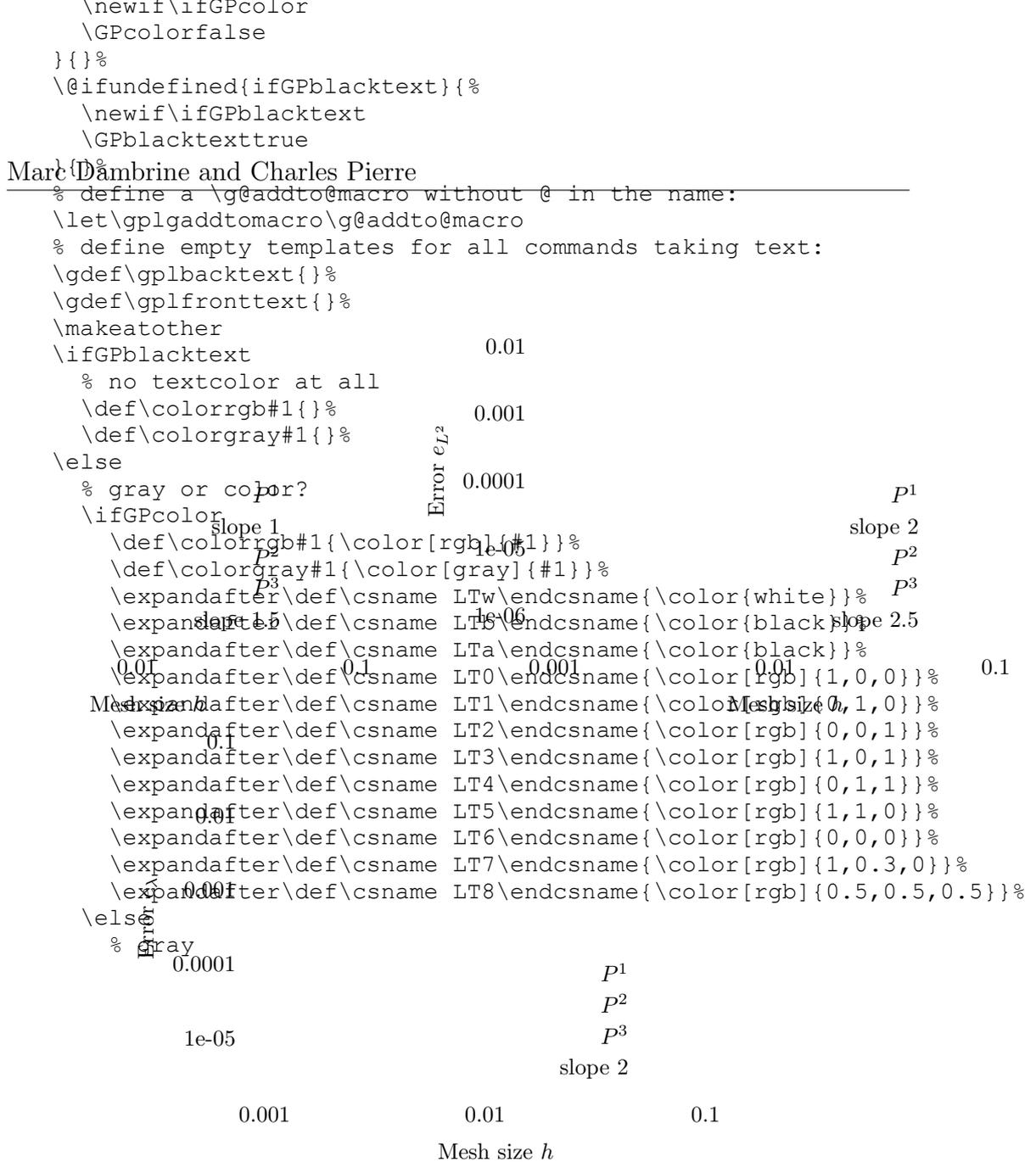}
  \end{tabular}\hspace{-2cm}
  \caption{2D case, convergence towards the 4$^th$ non zero eigenvalue on the disk} 
  \label{fig:ventcel-vp4-disk}
\end{figure}

\begin{figure}[h!]
  \centering
  \begin{tabular}{c}
    \hspace{-2.6cm}
    \input{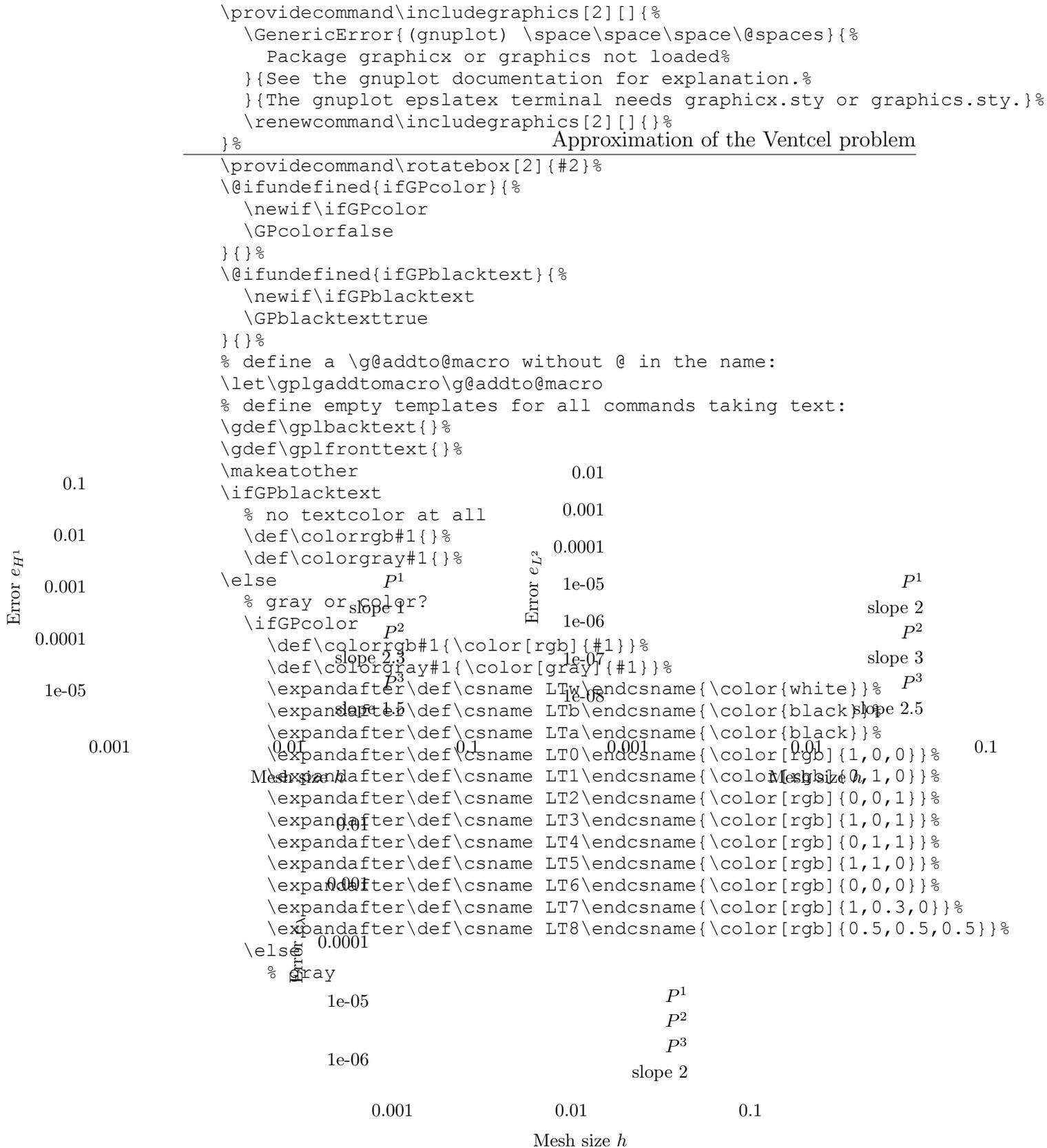}
  \end{tabular}\hspace{-2cm}
  \caption{Same as fihure \ref{fig:ventcel-vp4-disk} for the modified Ventcel problem}
\end{figure}

\bibliographystyle{abbrv}
\bibliography{biblio}

\def\cprime{$'$}
\begin{thebibliography}{1}

\bibitem{babuska-osborne-eigen}
I.~Babu{\v{s}}ka and J.~E. Osborn.
\newblock Estimates for the errors in eigenvalue and eigenvector approximation
  by {G}alerkin methods, with particular attention to the case of multiple
  eigenvalues.
\newblock {\em SIAM J. Numer. Anal.}, 24(6):1249--1276, 1987.

\bibitem{ciarlet-ef-elliptic}
P.~G. Ciarlet.
\newblock {\em The finite element method for elliptic problems}, volume~40 of
  {\em Classics in Applied Mathematics}.
\newblock Society for Industrial and Applied Mathematics (SIAM), Philadelphia,
  PA, 2002.
\newblock Reprint of the 1978 original [North-Holland, Amsterdam; MR0520174 (58
  \#25001)].

\bibitem{dambrine-kateb-2012}
M.~Dambrine and D.~Kateb.
\newblock Persistency of wellposedness of {V}entcel's boundary value problem
  under shape deformations.
\newblock {\em J. Math. Anal. Appl.}, 394(1):129--138, 2012.

\bibitem{dambrine-kateb-lamboley}
M.~Dambrine, D.~Kateb, and J.~Lamboley.
\newblock An extremal eigenvalue problem for the ventcel-laplace operator.
\newblock {\em ARXIV Preprint}, 2014.

\bibitem{Demlow}
A.~Demlow.
\newblock Higher-order finite element methods and pointwise error estimates for
  elliptic problems on surfaces.
\newblock {\em SIAM J. Numer. Anal.}, 47(2):805--827, 2009.

\bibitem{saad-vp}
Y.~Saad.
\newblock {\em Numerical methods for large eigenvalue problems}.
\newblock Algorithms and Architectures for Advanced Scientific Computing.
  Manchester University Press, Manchester, 1992.

\bibitem{saad}
Y.~Saad.
\newblock {\em Iterative methods for sparse linear systems}.
\newblock Society for Industrial and Applied Mathematics, Philadelphia, PA,
  second edition, 2003.

\bibitem{ventcel-56}
A.~D. Ventcel{\cprime}.
\newblock Semigroups of operators that correspond to a generalized differential
  operator of second order.
\newblock {\em Dokl. Akad. Nauk SSSR (N.S.)}, 111:269--272, 1956.

\bibitem{ventcel-59}
A.~D. Ventcel{\cprime}.
\newblock On boundary conditions for multi-dimensional diffusion processes.
\newblock {\em Theor. Probability Appl.}, 4:164--177, 1959.

\end{thebibliography}

\end{document}